\documentstyle[sprocl,amssymb]{article}

\newcommand{\A}{{\bf \cal A}}

\newcommand{\E}{{\Bbb E}}
\newcommand{\F}{{\cal F}}

\newcommand{\R}{{\Bbb R}}

\newcommand{\half}{{  {1\over 2}  }}

\newtheorem{theorem}{Theorem}[section]
\newtheorem{proposition}[theorem]{Proposition}
\newtheorem{lemma}[theorem]{Lemma}
\newtheorem{corollary}[theorem]{Corollary}

\def\exp{{\rm e}}

\begin{document}

\title{CONCERNING THE GEOMETRY OF STOCHASTIC DIFFERENTIAL EQUATIONS 
AND STOCHASTIC FLOWS}

\author{K.D.  ELWORTHY} 
\address{Department of Mathematics, University of Warwick, \\
Coventry CV4 7AL, UK}

\author{YVES  LE JAN }
\address{D\'epartment de Math\'ematique, Universit\'e Paris Sud,\\
 91405 Orsay, France}

\author{XUE-MEI  LI}
\address{Department of Mathematics, University of Warwick, \\
Coventry CV4 7AL, UK}

\maketitle

\begin{abstract}
Le Jan and Watanabe showed that a non-degenerate stochastic
flow $\{\xi_t: t\ge 0\}$ on a manifold $M$ determines a connection 
on $M$. This connection is characterized here and shown to be
 the Levi-Civita connection for gradient systems. This both explains
 why such systems have useful properties and allows us to extend these
properties to more general systems. Topics described here include:
moment estimates for $T\xi_t$, a Weitzenb\"ock formula for the generator
of the semigroup on p-forms induced by the flow,  a Bismut type
formula for $d\log p_t$ in terms of an arbitrary metric connection,
 and a generalized Bochner vanishing theorem.
\end{abstract}

\section{Introduction and Notations}

{\bf A.} Consider a Stratonovich stochastic differential equation
\begin{equation}\label{1}
dx_t=X(x_t)\circ dB_t+A(x_t)dt
\end{equation}
on an n-dimensional $C^\infty$ manifold $M$, e.g. $M=\R^n$.
Here $A$ is a $C^\infty$ vector field on $M$, so $A(x)$ lies in the
 tangent space $T_xM$ to $M$ at $x$ for each $x\in M$, while
 $X(x)\in {\cal L}(\R^m;T_xM)$, the space of linear maps of $\R^m$ to $T_xM$, 
for $x\in M$, and is $C^\infty$ in $x$. The noise $B_\cdot$ is a 
Brownian motion on $\R^m$ defined on a probability space $\{\Omega, \F, P\}$.

For each $e\in \R^m$ let $X^e$ be the vector field  given by 
$X^e(x)=X(x)(e)$. Recall that for each given $x_0\in M$ equation
(\ref{1}) has a maximal solution $\{\xi_t(x_0): 0\le t<\zeta(x_0)\}$,
defined up to an explosion time $\zeta(x_0)$, and unique up to
equivalence.  The  solutions form a Markov
process on $M$. Let $\{P_t^0: t\ge 0\}$ be the associated 
(sub)-Markovian semigroup, and let $\A$ be the infinitesimal generator. In
this article we shall assume that (\ref{1}) is {\it non-degenerate},
i.e. $X(x):\R^m\to T_xM$ is surjective for each $x$, or equivalently
that $\A$ is elliptic. Then a Riemannian metric is induced on
$M$ with inner product $<,>_x$ on $T_xM$ given by
$<X(x)e_1, X(x)e_2>_x=<e_1,e_2>_{\R^m}$ provided that $e_1$, $e_2$
are orthogonal to $N(x)$, the kernel of $X(x)$ in $\R^m$. The generator
has the form
\begin{equation}\label{2}
\A(f)(x)=\half \Delta^0f(x)
+\left<\half \sum_1^m \nabla X^{e_i}(X^{e_i}(x))+A(x), 
\hbox{ grad }f(x)\right>_x,
\end{equation}
where $e_1,\dots, e_m$ is an orthonormal basis for $\R^m$. Here
 $\nabla$ denotes covariant differentiation  with respect to 
the Levi-Civita connection,
so $\nabla X^{e_i}$ is a linear map of tangent vectors to tangent vectors,
$\nabla X^{e_i}(v)\equiv \nabla_v X^{e_i}$, and $\Delta^0$ is the 
Laplace-Beltrami operator on functions: 
$\Delta^0f =\hbox{trace }\nabla (\hbox{grad}f)$.

\bigskip

\noindent {\bf B.}
 Our motivating examples are the {\it gradient Brownian systems}.
Here we have an immersion: $g: M\to \R^m$, e.g. the inclusion of
the space of $S^n$ in $\R^{n+1}$ (with $m=n+1$), and  $X(x):\R^m\to T_xM$
is the orthogonal projection using $T_xg$ to identify $T_xM$ with a
subspace of $\R^m$. The Riemannian inner product $<,>_x$ is just that
which makes $T_xg$ an isometry. Set $Y(x)=T_xg: T_xM\to \R^m$. Let
$Z$ be a vector field then $Y(x)Z(x)\in \R^m$ for each $x$, giving
$Y(\cdot)Z(\cdot): M\to \R^m$, with differential
$d(Y(\cdot)Z(\cdot)): T_xM\to \R^m$, $x\in M$. It is a fundamental
result that if we project this differential  to $T_xM$ we obtain the
 Levi-Civita covariant derivative of $Z$ in the direction of $v$, i.e.
\begin{equation}\label{3}
\nabla Z(v)=X(x)\left[d(Y(\cdot)Z(\cdot))_x(v)\right], \hskip 18pt v\in T_xM,
\end{equation}
 e.g. see \cite{Kobayashi-NomizuII}.

\bigskip

 Consider the special case $Z(x)=X^e(x)$ some $e\in \R^m$. Then by 
(\ref{3}), for any $v\in T_xM$,
$$\nabla X^e(v)=X(x)\left[d\left(Y(\cdot)X(\cdot)e\right)(v)\right].$$
But  $Y(x)X(x)e=e-P_N(x)=P_T(x)$ say, where $P_N(x):\R^m\to \R^m$
is the orthogonal projection onto the normal space $N(x)$ at $x$,
and so, e.g. by differentiating the identity $P_T(x)e=P_T(x)P_T(x)e$,
we see that if $P_T(x)e=e$, i.e. if $e\in \hbox{Image } Y(x)$, then 
$\nabla X^e(v)=0$ all $v\in T_xM$ (for another proof see \S2A below).
Alternatively this can be seen from the fact that $\nabla X^\cdot(\cdot)$
is essentially the shape operator of the immersion. See e.g. \cite{ELbook}.
In particular from this we can conclude that the term 
$\sum_1^m\nabla X^{e_i}(x)(X^{e_i}(x))$ in (\ref{2}) vanishes so that
$\A f(x)=\half \Delta^0f(x)+ \left<A(x), \hbox{grad } f(x)\right>_x$.
These identities are behind the fact that  gradient systems have
particularly nice properties from the point of view of their solution flows,
see e.g. \cite{Kusuoka88}, \cite{EL-LI}, \cite{EL-RO94}, \cite{EL-YOR},
\cite{application}, and from the point of view of their It\^o maps 
\cite{Aida-EL95}.
Here we shall show that many of these constructions and properties are 
also true for general non-degenerate systems provided that we use connections
with torsion. Our starting point is:

\begin{theorem}\label{th:1}
For an arbitrary non-degenerate SDE (\ref{1}) there is a unique
affine connection $\tilde \nabla$ on $M$ such that
\begin{equation}\label{4}
\left(\tilde \nabla X^e\right)(v)=0, \hskip 18pt \hbox{all } v\in T_xM,
e\in [\hbox{ker}X(x)]^\perp.
\end{equation}
It is given by $\tilde \nabla Z(v)=\breve \nabla Z(v)$ for 
\begin{equation}\label{5}
\breve \nabla Z(v)=X(x)d[Y(\cdot)Z(\cdot)](v), \hskip 24pt v\in T_xM
\end{equation}
for $Y(x): T_xM \to \R^m$ the adjoint of $X(x)$, and is metric.
\end{theorem}
This is in fact the connection defined by LeJan and Watanabe 
\cite{LeJan-Watanabe82}.

\bigskip

\noindent{\bf C.} The scheme of the paper is as follows: Theorem 1 is proved
in \S 2 together with criteria for $\breve \nabla$ to be the Levi-Civita
connection and to be torsion-skew symmetric; in \S 3 we extend results
of \cite{EL-YOR} on the conditional expectation of the derivative flow
$T_{x_0}\xi_t$ given $\{\xi_t(x_0): 0\le t\le T\}$, i.e.
filtering out the extraneous noise; in \S 4 the 'spectral positivity'
estimates of \cite{application}, see also \cite{EL-RO94}, for moment
exponents are extended to S.D.E. with $\breve \nabla$ torsion skew
 symmetric; and in \S5 we give an expression for the generator
 of the semigroup $P_t^q$ on $q$-forms given by
 $P_t^q(\phi)=\E\xi_t^*(\phi)$ of the
form $P_t^q=-(\bar \delta d+d\bar\delta)$ and a Weitzenb\"ock formula.
An expression for the curvature is derived in Appendix I.

\bigskip

\noindent{\bf {Remark:}}

For simplicity in this expository article we mainly treat
 equations like (\ref{1}) with finite dimensional noise whereas
 stochastic flows correspond  canonically to 
Gaussian measures on the space of vector fields of $M$,
 \cite{Baxendale82}, \cite{Kunitabook}, \cite{LeJan-Watanabe82},
 which may have support on an infinite dimensional
space. Essentially this means that $\R^m$ should be replaced by a Hilbert
space of vector fields with $X(x)$ the evaluation map (the major role is 
then taken by the reproducing kernel of the Gaussian measure) 
\cite{Baxendale76}. 
See Appendix II. More generally
Gaussian measures on Hilbert spaces of sections of a vector bundle determines
a connection on that bundle (and all metric connections arise this way, 
see \S 2H below). Finally we also restrict ourselves here to non-degenerate
SDE, but a degenerate SDE induces in the same way a differential operator,
a 'semi-connection'. These aspects and other more detailed results will
be treated in a forthcoming article. See also \cite{EL-LJ-LI96}.

We are grateful to Profs. N. Ikeda and Z. Ma for helpful comments.
 For somewhat related work see \cite{Accardi-Mohari96}.

\section{Existence and basic properties}

\bigskip

\noindent {\it Proof of Theorem \ref{th:1}  }

Let $\breve \nabla$ be defined by (\ref{5}). It is easy to see that
it has the linearity and derivation properties, which ensures that it
is a connection.
 Let $\tilde \nabla$ be any affine connection on $M$, $Z$ a vector field,
and let $v\in T_{x_0}M$. Then 
\begin{equation}\label{7}
Z(x)=X(x)Y(x)Z(x), \hskip 18pt x\in M,
\end{equation}
whence $\tilde \nabla Z(v)=\tilde \nabla X(v)Y(x_0)Z(x_0)+\breve \nabla Z(v)$
 using (\ref{5}). 
 Setting $\tilde e=Y(x_0)Z(x_0)$ we see
\begin{equation}\label{8}
\tilde \nabla Z(v)=\tilde \nabla X^{\tilde e}(v)+\breve \nabla Z(v).
\end{equation}
Taking $\tilde \nabla=\breve \nabla$, since $Z(x_0)$ is arbitrary we see
$\breve \nabla$ satisfies the defining criterion (\ref{4}), giving
existence. Assuming $\tilde \nabla $ satisfies (\ref{4}) we see 
$\tilde \nabla Z(v)=\breve \nabla Z(v)$, giving uniqueness. To check
that $\breve \nabla$ is metric it is enough to show that
$$d\left(<Z(\cdot), Z(\cdot)>\right)(v)
=2 <\breve \nabla Z(v), Z(x_0)>_{x_0}.$$
In fact
\begin{eqnarray*}
<\breve \nabla Z(v),\, Z(x_0)>_{x_0}&=&\left<
d[Y(\cdot)Z(\cdot)](v),\, Y(x_0)Z(x_0)\right>_{\R^m}\\
&=&{1\over 2} d\left<Y(\cdot)Z(\cdot),\, Y(\cdot)Z(\cdot)\right>_{\R^m}(v)\\
&=&\half d\left<Z(\cdot),\, Z(\cdot)\right>_\cdot (v).
\end{eqnarray*}
\hfill //

\bigskip

\noindent{\bf Remark:}

 Note that $Y(x)Z(x)=\sum_1^m<X^{e_i}(x),Z(x)>e_i $ and by (\ref{7}) and the
 equation which follows:
\begin{equation}\label{201}
\begin{array}{ll}
\breve \nabla Z(v)&=\sum_1^m X^{e_i}d<X^{e_i}, Z>(v)\\
&=\tilde\nabla Z(v)-\tilde \nabla_v X^{e_i}<X^{e_i}(x_0), Z(x_0)>
\end{array}\end{equation}
for any affine connection $\tilde \nabla$ on $M$.
\bigskip

\noindent{\bf B.}
 In a local chart about $x_0\in M$ we can take $\tilde \nabla$ in 
the above proof to be the usual derivative so that (\ref{8})
becomes
$$DZ(x_0)(v)=DX^{\tilde e}(x_0)(v)+\breve \nabla Z(v)$$
where $\tilde e=Y(x_0)Z(x_0)$, 
(using local representations for $Z$, $X^{\tilde e}$, and $v$).
But for $\breve \Gamma$ the Christoffel symbol of $\breve \nabla $ in
 our chart
$$\breve \nabla Z(v)=DZ(x_0)(v)+\breve \Gamma(x_0)(v, Z(x_0))$$
giving
\begin{equation}\label {9}
\breve \Gamma (x_0)(v,w)
=-DX(x_0)(v)\left(Y(x_0)w\right), \hskip 15pt v,w\in R^n.
\end{equation}
Equivalently 
\begin{equation}
\breve \Gamma_{jk}^i=-\sum_{r=1}^m\sum_{l=1}^n
 {\partial X(x_0)^{r,i}\over\partial x^j} X(x_0)^{r,\ell}g_{k\ell},
\end{equation}
where $\left\{ X(x)^{r,i} \right\}$, $\{1\le i\le n\}$, $\{1\le r\le m\}$
is the matrix representing $X(x): \R^m\to \R$, i.e.
$X(x)^{r,i}=<X(e_r),  f_i>$ for $\{e_i\}$ and $\{f_i\}$
orthonormal bases for $\R^m$ and $T_xM$ respectively, and
$\{g_{k\ell}\}$ the metric tensor. This shows that $\breve \nabla$
is the LeJan-Watanabe connection defined in \cite{LeJan-Watanabe82}.

\bigskip

\noindent {\bf C.} Equivalent definitions and properties.

\begin{lemma}\label{le:2.1} For any orthonormal base $\{e_i\}$ of $\R^m$ and 
$v\in T_{x_0}M$ we have
\begin{equation}\label{10}
(i) \hskip 28pt 
\breve \nabla Z(v)={d\over dt}
\sum_1^m X^{e_i}(x_0)\left.\left
< Z(\sigma(t)), X^{e_i}(\sigma(t))\right>_{\sigma(t)}\right|_{t=0}
\end{equation}
where $\sigma: [-\delta, \delta]\to M$ is a $C^1$ curve with $\sigma(0)=x_0$
and $\dot \sigma(0)=v$.

$$ (ii)\hskip 30pt
\breve \nabla Z(v)=\sum_1^m [X^i, V](x_0) <X^i(x_0), Z(x_0)>+[V, Z](x_0)$$
where $V$ is any smooth vector field with $V(x_0)=v$.
\end{lemma}

\noindent{\bf Proof.}
Since $\breve \nabla $ is metric the right hand side of (\ref{10})
is just
$$\sum_1^m X^{e_i}(x_0)\left\{
\left<\breve \nabla Z(v),\, X^{e_i}(x_0)\right>_{x_0}
+\left<Z(x_0),\, \breve\nabla X^{e_i}(v)\right>_{x_0}\right\}.$$
This is independent of the choice of basis. Choose $\{e_i\}$
so that $e_1\dots, e_n$ span $[\hbox{ker}X(x_0)]^\perp$, i.e.
are in the image of $Y(x_0)$. Then $X^{e_i}(x_0)=0$ if $i>n$
while $\breve\nabla X^{e_i}(v)=0$ if $1\le i\le n$ by definition of
$\breve \nabla$. Since $X^{e_i}(x_0), 1\le i\le n$, form an 
orthonormal base for $T_{x_0}M$ the result (i) follows. 

For (ii) write
$$[V, Z]= [ V,\,\,  \sum_1^m <X^i, Z> X^i \,]$$
and expand. The use of (\ref{201}) yields (ii).

\hfill  //

\bigskip
By a similar proof to that above, we obtain a necessary and sufficient 
condition for a connection to be a metric connection: for simplicity
write $X^i\equiv X^{e_i}$,

\begin{lemma} 
A connection $\tilde \nabla$ is a metric connection if and only if
\begin{equation}\label{11}
\sum_1^m X^{e_i}<Z, \tilde\nabla_v X^{e_i}>+
\sum_1^m \tilde\nabla_v X^{e_i} <Z, X^{e_i}>=0,
\end{equation}
for all vector fields $Z$.
\end{lemma}

\noindent{\bf Proof.} 
Take $v\in T_xM$ . If $\tilde \nabla $ is metric then
\begin{eqnarray*}
&&d<Z(\cdot), Z(\cdot)>(v)=\sum_1^m d\left(<Z, X^{i}><Z, X^{i}>\right)(v)\\
=&&2 <\tilde \nabla Z(v), Z>+2\sum_1^m <Z, \tilde\nabla X^{i}(v)>
 <Z, X^{i}>\\
\end{eqnarray*}
giving (\ref{11}) by polarization.  
Now suppose  (\ref{11}) holds for a connection $\tilde \nabla$, then
\begin{equation}\label{112}
\sum_1^m<Z, X^i><Z,\tilde \nabla X^i(v)>=0.
\end{equation}
On the other hand, by (\ref{7})
\begin{eqnarray*}
&&\tilde \nabla Z(v)= \tilde \nabla_v Y(x)Z(x)+X(x) d[Y(x)Z(x)](v)\\
=&& \sum_1^m \tilde\nabla_v X^{i} <Z, X^{i}>+
   \sum_1^m X^{i} d < Z(-), X^{i}(-)>(v)\\
\end{eqnarray*}
giving
$$\sum_1^m X^{i} d<Z(-), X^{i}(-)>(v)=
\tilde \nabla Z(v) -\sum_1^m \tilde \nabla_v X^{i} <Z, X^{i}>. $$
consequently
\begin{eqnarray*}
&&d<Z(\cdot), Z(\cdot)>(v)= \sum_1^m d< Z(-), X^i(-)>^2(v)\\
=&& 2\sum_1^m <Z(x), X^i(x)> d< Z(-), X^i(-)>(v)\\
=&&2<Z, \tilde \nabla Z(v)> -2\sum_1^m <Z, \tilde\nabla_v X^{i}><Z, X^{i}>\\
=&&2  <Z, \tilde \nabla_v Z>,
\end{eqnarray*}
using (\ref{112}), and so $\tilde \nabla $ is 
a metric connection. \hfill //

\bigskip

\noindent {\bf D.} 
Recall that for any connection $\tilde \nabla$ on $M$ the
 torsion is a bilinear map from tangent vectors to tangent vectors,
$\tilde T: TM\oplus TM\to TM$, given by
\begin{equation}\label{12}
\tilde T(U(x_0),V(x_0))=\tilde \nabla V(U(x_0))-\tilde\nabla U(V(x_0))
-[U,V](x_0)
\end{equation}
for vector fields $U$, $V$.

Let $v_1, v_2\in T_{x_0}M$. There are the vector fields $Z^{v_1}$,
$Z^{v_2}$ given by
$$Z^{v_i}=X(x)Y(x_0)v_i, \hskip 18pt  i=1,2.$$
By definition
$$\breve \nabla Z^{v_i}(v)=0,   \hskip 38pt  \hbox{any } v\in T_{x_0}M.$$
Thus
\begin{equation}\label{13}
\breve T(v_1, v_2)=-\left[Z^{v_1},\, Z^{v_2}\right](x_0).
\end{equation}

Alternatively using the Levi-Civita connection in  (\ref{5})
\begin{equation}\label{14}
\breve\nabla Z(v)=X(x_0)\nabla Y(v) Z(x_0)+\nabla Z(v)
\end{equation}
whence by (\ref{12}) 
\begin{eqnarray*}
\breve T(v_1,v_2)&=&X(x_0)\left(\nabla Y(v_1)(v_2)-\nabla Y(v_2)(v_1)\right)\\
&&+\nabla Z^{v_2}(v_1)-\nabla Z^{v_1}(v_2)-[Z^{v_1}, Z^{v_2}].
\end{eqnarray*}
Thus by (\ref{12}) and the standard formula for exterior differentiation:
\begin{equation}\label{15}
\breve T(v_1,v_2)=X(x_0)dY(v_1,v_2), \hskip 20pt v_1, v_2\in T_{x_0}M.
\end{equation}

\bigskip

\noindent{\bf E.}
For any connection $\tilde\nabla $ on $M$, there is an adjoint
connection $\tilde\nabla^\prime$ on $M$ defined by
\begin{eqnarray*}
\tilde\nabla^\prime Z(v)&=&\tilde \nabla Z(v)-\tilde T(v,\,Z(x_0))\\
&=&\tilde \nabla V(Z(x_0))- [Z, V](x_0).
\end{eqnarray*}
Here $V$ is a vector field such that $V(x_0)=v$.  In terms of
 Christoffel symbols (\cite{Driver92}) this is equivalent to
 $ \tilde{ \Gamma^\prime}_{jk}^i=\tilde \Gamma_{kj}^i$. 
If $\hat \nabla$ denotes adjoint of $\breve \nabla$ we see
that $\hat \nabla Z(v)=[ Z^v, Z](x_0)$.

A connection $\tilde \nabla$ on a Riemannian manifold $M$ is
{\it torsion skew symmetric}, see \cite{Driver92}, if 
$u\to \tilde T(u,v)$ is skew symmetric as a map
$T_{x_0}M\to T_{x_0}M$ for all $v\in T_{x_0}M$, all $x_0\in M$.
 We have:

\bigskip

\begin{lemma}\label{le:1E}
A metric connection $\tilde \nabla$ on a Riemannian manifold $M$ is
torsion skew symmetric if and only if its adjoint connection is metric. 
If so the geodesics for $\tilde \nabla$ are those of the Levi-Civita
connection and the (usual) Laplace-Beltrami operator acting on a 
function $f$, $\Delta^0f$, is given by the trace of 
$\tilde\nabla(\hbox{grad} f)$.
\end{lemma}

\noindent{\bf Proof.} See \cite{Driver92} and also \cite{Kobayashi-NomizuI}
 (the last part also comes from the next proposition). 

\bigskip

\begin{proposition}\label{pr:1E}
The connection $\breve \nabla$ is
\begin{enumerate}
\item
the Levi-Civita connection if and only if $\nabla Z^v$ vanishes at $x_0$
for all $v\in T_{x_0}M$.
\item
torsion skew symmetric if and only if 
$ {\left. \nabla Z^v \right|}_{T_{x_0}M}: T_{x_0}M\to T_{x_0}M$
is skew symmetric, all $v\in T_{x_0}M$, or equivalently
$\nabla_v Z^w+\nabla_wZ^v=0$  for any $w,v\in T_xM$,
or  $\breve\nabla _UV+\breve\nabla_VU=\nabla_UV+\nabla_VU$
for all  vector fields $U$ and $V$.
\end{enumerate}
Also it is Levi-Civita if and only if $X(x)dY(u,v)=0$ for all
$u,v\in T_xM$, all $x\in M$.
\end{proposition}

\noindent{\bf Proof.}
The first part comes from the defining property of $\breve \nabla$ and
the third part comes from (\ref{15}).  For the second part, first observe by
  the definition  of torsion
$$\breve T(u,v)=\breve \nabla_v Z^u- \nabla_v Z^u
-\left[\breve \nabla_u Z^v-\nabla_u Z^v\right].$$
and so by (\ref{8}):

\begin{equation}
\breve T(u,v)=\sum_1^m X^{i}<v, \nabla X^{i}(u)>
-\sum_1^m X^{i}<u, \nabla X^{i}(v)>.
\end{equation}
 We have:
\begin{eqnarray*}
&&<\breve T(u,v), w>\\
=&&\sum_1^m <X^i, w><v, \nabla X^i(u)>
-\sum_1^m<X^i,w><u, \nabla X^i(v)>.
\end{eqnarray*}
However the second term is anti-symmetric in $u$ and $w$ by (\ref{11}).
Thus

\begin{eqnarray*}
&&<\breve T(u,v), w>+<\breve T(w,v), u>\\
=&& \sum_1^m <X^i, w><v,\nabla X^i(u)>+\sum_1^m <X^i, u><v,\nabla X^i(w)>\\
=&& \left<\nabla_u Z^w, v\right>+ \left<\nabla_w Z^u, v\right>\\
=&&- <w, \nabla Z^v(u)>-<u, \nabla Z^v(w)>,
\end{eqnarray*}
since $d<Z^w, Z^v>(u)=0$ and $d<Z^u, Z^v>(w)=0$. 

Also if $U$ and $V$ are vector fields, by (\ref{201})
$$\breve\nabla_VU=\sum_1^m X^i <U, \nabla_V X^i>+\nabla_VU$$
and so 
$$\breve\nabla_VU+\breve\nabla_UV=\nabla_VU+\nabla_UV+A$$
for $$A=\sum_1^m X^i<U, \nabla_V X^i>+ \sum_1^m X^i<V, \nabla_U X^i>.$$
But $\breve T$ is skew symmetric if and only if $A\equiv 0$.
\hfill //

\bigskip

\begin{corollary}\label{2.5}
If $\breve \nabla$ is torsion skew symmetric then  
 $$\breve T(u,v)= 2\sum_{i=1}^m X^{i}<u, \nabla X^{i}(v)>$$ and the
Levi -Civita connection can be expressed in terms of the  LeJan-Watanabe
 connection by: 
\begin{equation}\label{117}
\nabla  Z(v)=\breve \nabla Z(v)-\half \breve T(Z(x_0),v).
\end{equation}
In particular $\nabla X^i(X^i)=0$ for each $i$.
\end{corollary}

\bigskip

\noindent
{\bf Remark:} Most of the results for gradient Brownian systems carry
over to the case when $\breve \nabla$ is torsion free and,
with some adaptation, to the torsion skew symmetric case or even
more generally.

\bigskip

\noindent{\bf F.}
Let $f: M\to \R$ be  $C^2$.  Then It\^o's formula gives
\begin{eqnarray*}
f(\xi_t(x))&=&f(x_0)+\int_0^t df\left(X(\xi_s(x_0))dB_s\right)\\
&&+\half \int_0^t \hbox{trace} \breve \nabla (\hbox{grad} f)(\xi_s(x_0))ds\\ 
&&+\int_0^t A(\xi_s(x_0))ds, \hskip 15pt 0\le t< \zeta(x_0)
\end{eqnarray*}
since the Stratonovich term
 $\sum_1^m \breve \nabla X^{e_i}(X^{e_i}(x))$
vanishes. Thus as shown in \cite{LeJan-Watanabe82}, the generator
is given by
\begin{equation}\label{118}
\A^0(f)=\half \hbox{trace} \breve \nabla (\hbox{grad} f)+ <A(\cdot),
 \hbox{grad}f>.
\end{equation}
Note also that the vanishing of the Stratonovich term means that
{\it (\ref{1}) can be considered as an It\^o equation w.r.t.}
 $\breve \nabla$, e.g. see \cite{ELbook}, and {\it the solutions 
$\{\xi_t(x_0): t\ge 0\}$ will be $\breve \nabla$-martingales if} $A\equiv 0$, 
\cite{Emerybook}. Furthermore by  Corollary \ref{2.5}
 if $\breve \nabla$  is torsion skew symmetric 
(\ref{1}) will be an It\^o equation for  the Levi-Civita connection 
and the solution will be a Brownian  motion with drift $A$.

\bigskip

\noindent {\bf G.}
 Example: Invariant SDE on Lie groups: c.f. \cite{Driver92}.
Let $M$ be a Lie group and suppose (\ref{1}) is a left invariant SDE, 
with $A=0$ for simplicity. For $g\in G$ let $R_g: G\to G$ and
$L_g: G\to G$ be right and left translations by $g$. Then
$$L_gX(x)(e)=X(gx)e, \hskip 24pt g,x\in G, e\in \R^m.$$
We can suppose $m=n$ since $\hbox{Ker} X(x)$ is independent of $x$.
The metric induced on $G$ will be left invariant. We can treat
$X(id): \R^m\to T_{id}G$ as an identification of $\R^m$ with
the Lie algebra ${\cal G} =T_{id} G$ of $G$, and then $Y$ becomes
the Maurer-Cartan form. For $v\in T_{x_0}G$ the vector field 
$Z^v=X(\cdot)Y(x_0)(v)$ of \S 1D is just the left-invariant vector field
through $v$. If $\tilde \nabla$ is the flat left
invariant connection on $G$ then $\tilde \nabla Z^v\equiv 0$,
 and so by definition $\breve \nabla=\tilde \nabla$. The torsion
\begin{eqnarray*}
\breve T(v_1,v_2)&=&-\left[Z^{v_1},\, Z^{v_2}\right](x_0)\\
&=&X(x)dY(v_1,v_2)
\end{eqnarray*}
by (\ref{13}) and (\ref{15}).

Recall that for $\alpha\in {\cal G}$,
$$ad(\alpha): {\cal G}\to {\cal G}$$
is given by
$$ ad(\alpha)\beta=[\alpha, \beta].$$
Taking $x_0=id\in G$ we see $\breve T(v_1,v_2)=-ad(v_1)(v_2)$
and so $\breve \nabla $ is torsion skew symmetric if and only 
if $ad(v_1)$ is skew symmetric for all $v_1\in {\cal G}$. From
Lemma 7.2 of \cite{Milnor76} we know this holds if and only
if the metric on $G$ is bi-invariant (which is only possible
if $G$ is isomorphic to the product of compact group and a 
commutative group). Indeed from the proof of Lemma 7.2 and
7.1 of \cite{Milnor76} we see that the adjoint connection
is the flat right invariant connection, which is a metric
connection for the right invariant metric

$$<v_1,v_2>^\prime_{x_0}\equiv 
\left< TR_{x_0}^{-1}(v_1),  TR_{x_0}^{-1}(v_2)\right>_{id}.$$

\bigskip

\noindent{\bf H.}
There is a natural correspondence between S.D.E.'s (\ref{1}) with
$A\equiv 0$ and smooth maps of $M$ into the Grassmanian of n-planes
in $\R^m$ classifying $TM$. The connection $\breve \nabla$ is
the pull back of the universal connection on the Stiefel bundle
over $M$, described in \cite{Narasimhan-Ramanan61}. From there it follows
that {\it every metric connection on $M$ can be obtained as 
$\breve \nabla$ for some} S.D.E. (\ref{1}), see \cite{EL-LJ-LI}.

\bigskip

For a diffusion on $M$, with $n= \hbox{dim } M>1$,
 with generator $\half \Delta+L_Z$, for some smooth
vector field $Z$, Ikeda and Watanabe showed how to construct a
metric connection $\tilde \nabla$ on $M$ such that the diffusion 
process (from any point $x_0$ of $M$), is a $\tilde \nabla$-martingale
(it is the stochastic development of an n-dimensional Brownian motion).
See \cite{Ikeda-Watanabe}. This $\tilde \nabla$ is not uniquely
determined. By the remark above $\tilde \nabla =\breve \nabla$
for some S.D.E. $dx_t=X(x_t)\circ dB_t$, again not uniquely determined. 
For this S.D.E. the generator satisfies
$\sum_{i=1}^m L_{X^i}L_{X^i}=\half \Delta +L_Z$. As T. Lyons has
pointed out to us such a construction is not in general possible when 
$\hbox{dim } M=1$.

\bigskip

\noindent{\bf I.}
We summarize here some of the notation being used:
\begin{eqnarray*}
N(x)&=&\hbox{Ker} X(x), \\
 Y(x)&=& X(x)^*: T_xM \to \R^m,\\
 Z^v&=&X(\cdot)Y(x_0)v,\hskip 6pt v\in T_{x_0}M\\
\nabla,&&\hbox{Levi-Civita connection, $R$, Ric, 
its curvature and Ricci}\\
&& \hbox{curvature;}\\
\tilde\nabla,&& \hbox{any connection, $\tilde R$, $\tilde Ric$, 
its curvature and Ricci }\\
&& \hbox{curvature,
 $\tilde{{\rm Ric}}^\#(v)=\sum_1^m \tilde{{\rm Ric}} (v, X^i(x))X^i(x)$},\\
&&\hbox{  and $\tilde T$ its torsion tensor}\\
\breve\nabla, && \hbox{LeJan-Watanabe connection, $\breve R$, $\breve Ric$, 
its curvature and}\\
&& \hbox{ Ricci curvature, and $\breve T$ its torsion tensor}\\
\hat\nabla, && \hbox{the adjoint  connection of $\breve \nabla$, 
$\hat R$, $\hat Ric$, its curvature and},\\
&&\hbox{ Ricci curvature, and $\hat T$ its torsion tensor}.\\
\end{eqnarray*}

\section{The Derivative flow}

{\bf A.} Let $N=\cup_{x\in M} N(x)$. It forms a Riemannian vector bundle
over $M$, (the normal bundle in the gradient case). Take any metric
connection on it, with parallel translation along a curve
 $\{\sigma(s): 0\le s\le t\}$ denoted by 
$\tilde{//_s}: N(\sigma(0))\to N(\sigma(s))$. Let
 $\breve{//_t}$ be parallel translation for $\breve \nabla$. Using $Y$ this
induces a parallel translation operator
 $$//_t
=Y(\sigma(t))\breve{//_t} X(\sigma(0)): 
N(\sigma(0))^\perp\to N(\sigma(t))^\perp,$$
which combines with $\tilde {//}_t$ on $N$ to give a parallel translation
in $M\times \R^m$, again written $\tilde {//}_t$, as an isometry
$$\tilde{//}_t: \R^m\to \R^m$$
depending on $\sigma$. Following \cite{EL-YOR}, set
\begin{equation}
\breve B_t:= \int_0^t \breve{//}_s^{-1} X(x_s)dB_s
\end{equation}
and 
\begin{equation}
\beta_t:=\int_0^t  \tilde{//_s}^{-1} K(x_s) dB_s
\end{equation}
where $K(x): \R^m\to \R^m$ is the orthogonal projection onto
 $N(x)$ and
$x_s=\xi_s(x_0)$. Finally set
$$\tilde B_t=Y(x_0)\breve B_t=\int_0^t \tilde{//_s}^{-1}Y(x_s)X(x_s)dB_s,$$
and  $\bar B_t=\tilde B_t+\beta_t$.

For any process $\{y_s: 0\le s< \zeta\}$ let
$\F^{y_\cdot}=\sigma\{y_s: 0\le s<\zeta\}$, but write 
$\F^{\xi_\cdot(x_0)}$ as $\F^{x_0}$. The following decomposition
 theorem is a direct analogue of the corresponding results in
\cite{EL-YOR} with the same proof:

\begin{theorem}\label{th:2A}

\begin{enumerate}
\item
$\F^{\breve B_\cdot}=\F^{x_0}$,
\item
$\{\bar B_t: 0\le t< \zeta\}$ is a Brownian motion on $\R^m$
with $B_t=\int_0^t \tilde{//_s}d\bar B_s$.
\end{enumerate}
In particular $\{\beta_t: 0\le t<\zeta\}$, when conditioned on
$\{\breve B_t: 0\le t<\zeta\}$  is a Brownian motion killed
at time $\zeta$ (so when $\zeta=\infty$, $\beta_\cdot$ and $\breve B_\cdot$
are independent Brownian motions).
\end{theorem}

\bigskip

\noindent{\bf B.}
 The derivative flow $T\xi_t$ on $TM$ is given by the covariant
equation
\begin{equation}\label{500}
\tilde D v_t=\tilde \nabla X(v_t)\circ dB_t+ \tilde \nabla A(v_t)dt-
\tilde T(v_t, X(x_t)\circ dB_t+A(x_t)dt)
\end{equation}
for $v_t=T\xi_t(v_0)$, along the paths of $\{\xi_t: 0\le t<\zeta\}$,
since for a $C^1$ map
 $\sigma: (-\delta, \delta)\times (-\delta, \delta) \to M$
$${\tilde D\over \partial s}{\partial\over\partial t } \sigma(s,t)=
{\tilde D\over \partial t}{\partial\over \partial s} \sigma(s,t)
+\tilde T( {\partial \sigma\over \partial s}, 
{\partial \sigma\over\partial t})$$
(e.g. see \cite{Milnorbook}).  Such covariant equations are
 described in \cite{ELbook}. 

Taking $\tilde \nabla$ to be the adjoint connection  $\hat \nabla$, since
\begin{equation}
\hat \nabla_U V=\breve \nabla _U V - \breve  T(U, V),
\end{equation}
we see 
\begin{equation}
\hat Dv_t=\breve \nabla X(v_t)\circ dB_t+\breve \nabla A(v_t)dt.
\end{equation}
To rewrite this as an It\^o equation (which means apply $\hat{//_t}^{-1}$ 
to both sides and consider the resulting It\^o equation in $T_{x_0}M$),
the correction term is
\begin{eqnarray*}
&&\half \sum_1^m \left[\breve \nabla X^i
\left(\breve \nabla X^i(v_t)\right)dt    +\hat \nabla_{X^i}
 \left(\breve\nabla X^i\right) \left(v_t\right)dt\right] \\
&=&\half \sum_1^m \left[\breve \nabla X^i
\left(\breve \nabla X^i(v_t)\right)dt +\breve \nabla^2 
 X^i\left(X^i, v_t\right)dt+\breve T
\left(\breve\nabla X^i(v_t), X^i\right)\right]\\
&=&\half \sum_{i=1}^m \left[\breve \nabla 
\left(\breve\nabla  X^i\left( X^i(\cdot)\right)\right)(v_t)
+\breve\nabla^2  X^i\left(X^i, v_t\right)-
\breve\nabla^2  X^i\left(v_t, X^i\right)\right]dt\\
&=&\half \sum_{i=1}^m \left[\breve \nabla 
\left(\breve\nabla  X^i\left( X^i(\cdot)\right)\right)(v_t)dt
-\half \breve{\hbox{Ric}}^\#(v_t)dt\right]
\end{eqnarray*}
as in \cite{ELflour}, \cite{EL-YOR}, where 
$\breve{\hbox{ Ric}}^\#(v)
=\sum_1^m \breve{\hbox{ R}}\left(v, X^i(x)\right)X^i(x)$
so that
$<\breve{\hbox{Ric}}^{\#}(v_1),v_2 >_x=\breve{\hbox{Ric}}(v_1,v_2)$
 for $v_1, v_2\in T_xM$. The first term vanishes as we saw in \S 1E from
 the definition of $\breve\nabla$. Thus
\begin{equation}\label{19}
\hat Dv_t=\breve \nabla X(v_t)dB_t-\half\breve {\hbox{Ric}}^\#(v_t)dt+
\breve \nabla A(v_t)dt.
\end{equation}

\bigskip

\noindent {\bf C.}
 We can now extend one of the main results of \cite{EL-YOR}.
If $\{u_t: 0\le t<\zeta\}$ is any process along $\{\xi_t: 0\le t<\zeta\}$
by
$\E\{u_t\,\chi_{t<\zeta(x_0)}|\F^{x_0}\}$ we mean
 $\tilde{//_t}\E\{\tilde{//_t}^{-1} u_t\, 
\chi_{t<\zeta(x_0)}|\F^{x_0}\}$.
 As pointed out by 
M. Emery this is independent of the connection $\tilde\nabla$ used
 to define $\tilde{//_t}$.

\begin{theorem}\label{th:2C}
Assume $|v_t|$ is integrable for each $t\ge 0$.   Set \newline
$v_t^{x_0}= \E\{T\xi_t(v_0)\chi_{t<\zeta(x_0)}|\F^{x_0}\}$. Then
$\{v_t^{x_0}\}$ satisfies the covariant equation
\begin{equation}\label{20}
\hat Dv_t^{x_0}=-\half\breve{\hbox{Ric}}^\#(v_t^{x_0})dt+
\breve \nabla A(v_t^{x_0})dt.
\end{equation}
along $\{x_t\}$ on $t<\zeta$.
\end{theorem}

\noindent
{\bf Proof.} First assume non-explosion.
Using Theorem \ref{th:2A} and rewriting (\ref{19}) as
\begin{eqnarray*}
\hat D v_t&=&\breve \nabla X(v_t)\tilde{//_t} d\bar B_t-
\half \breve{\hbox{Ric}}^\#(v_t)dt+\breve \nabla A(v_t)dt\\
&=&\breve \nabla X(v_t)\tilde{//_t} d\beta_t-
\half \breve{\hbox{Ric}}^\#(v_t)dt+\breve \nabla A(v_t)dt.\\
\end{eqnarray*}
The last step used the fact that 
 $$\breve \nabla X(v_t)\left(\tilde{//_t}d\tilde B_t\right)
=\breve \nabla X(v_t) \left(Y(x_t)X(x_t)dB_t\right)=0,$$
by definition of $\breve\nabla$.
But by theorem \ref{th:2A}
$$\E\left\{\int_0^t \breve \nabla X(v_s)\tilde{//_s}d \beta_s\,
|\, \F^{x_0}\right\}=0$$
and the result follows by the linearity and $\F^{x_0}$-measurability
of $\hbox{Ric}^{\#}_{x_t}$ and $\breve \nabla A_{x_t}$.
 If $\zeta(x_0)<\infty$,  let $\tau_D$ be the first exit time of
 $\xi(x_0)$ from a domain $D$ of $M$ with $D$ compact. The above argument
show (\ref{20}) holds on $t<\tau_D$. Now choose $D^i$ with
 $\tau_{D^i}\to \zeta$.
\hfill //

\bigskip

\noindent{\bf {Remark 2. }} 
 The integrability of $|v_t|$ is needed in
order for $v_t^{x_0}$ to be defined. It holds if $M$ is compact or
with  conditions on the growth of $|\nabla X|$, $|\nabla^2 X|$,
 and $|\nabla A|$
\cite{flow}, and is close to implying non-explosion of $\{x_t: t\ge 0\}$,
\cite{application}.

\bigskip

As an illustrative application there is the following extension
 of Bochner's vanishing theorem (however see Proposition \ref{4.3}
below):

\begin{corollary}\label{co:2C1}
Suppose $M$ is compact. If $M$ admits a vector field $A$ and a
metric  connection $\tilde\nabla$ whose adjoint 
connection preserves a  metric $<-,->^\prime$  such that 
$$\left<\tilde{\hbox{Ric}}^\#(v), v)\right>^\prime  
 >  2\left<\tilde\nabla A(v), v\right>^\prime
\hskip 18pt \hbox{all }  v\in TM, v\not =0,$$
then the cohomology group $H^1(M;\R)$ vanishes.
\end{corollary}

\noindent{\bf Proof.} 
Let $\phi$ be a closed smooth 1-form and $\sigma$ a singular
 1-cycle in $M$. By DeRham's theorem it is enough to show
 $\int_\sigma\phi=0$. According to \S2 H we can find an SDE (\ref{1})
with $\tilde \nabla =\breve \nabla$. Since $M$ is compact (\ref{1})
has a smooth solution flow $\{\xi_t: t\ge 0\}$ of diffeomorphisms of
 $M$. Then, by the continuity in $(t,x)\in \R(\ge 0)\times M$ of $\xi$,

$$\int_\sigma \phi dx =\int_{\xi_t\sigma}\phi dx
=\int_\sigma \xi_t^* \phi dx$$
Treating the case when $\sigma: [a,b]\to M$ this gives
\begin{eqnarray*}
\int_\sigma \phi dx &=&\E \int_a^b \phi_{\xi_t(\sigma(\theta))}
\left(T\xi_t(\dot \sigma(\theta))\right)d\theta\\
&=&\E \int_a^b \phi_{\xi_t(\sigma(\theta))}
\E \left\{T\xi_t(\dot \sigma(\theta)) |\F_t^{\sigma(\theta)}\right\}
d\theta\\
&=& \int_a^b \E\phi_{\xi_t(\sigma(\theta))}
 \breve W_t^A (\dot \sigma(\theta))d\theta
\end{eqnarray*}
where 
\[\left\{
\begin{array}{lll}
   {\hat D\over \partial t} \breve W_t^A(v_0)&=&-\half
\breve{\hbox{Ric}}^\# (\breve{W}_t^A(v_0))
+\breve \nabla A(\breve{W}_t^A(v_0))\\
\breve {W}_0^A(v_0)&=&v_0\in TM. 
\end{array}\right.  \]
by Theorem \ref{th:2C}. Thus 
$$|\int_\sigma \phi | \le 
\sup_x|\phi_x|^\prime \int_a^b |\breve {W}_t^A(\dot \sigma(\theta)|^\prime
d\theta.$$
However, for $v_0\in T_{x_0}M$,
$${d\over dt}|\breve W_t^A(v_0)|^{\prime, 2}
=2 \left<{\hat D\over \partial t}\breve W_t^A(v_0),
 \breve W_t^A(v_0)\right>^\prime.$$
So our assumptions imply that $\breve W_t^A(v_0)$ decays exponentially
as $t\to \infty$,  uniformly in $x_0\in M$, $v\in T_{x_0}M$ with 
$|v_0|^\prime=1$. Thus, letting $t\to \infty$, we see $\int_\sigma \phi=0$.
\hfill //

\bigskip

Next we give a version of Bismut's formula in this context, c.f.
 \cite{Driver92}.

\begin{corollary}\label{co:2C2}
 Let $\tilde \nabla$ be a metric connection
for a compact Riemannian manifold $M$.
 Let $p_t(x,y)$ be the  fundamental solution to
$${\partial u_t\over \partial t}
=\half \hbox{trace} \tilde \nabla (\hbox{grad } u_t)+L_A u_t.$$
Then 
$$d\log p_t(\cdot, y)(v_0)={1\over t}\E
\left\{\int_0^t
 \left<\tilde W_s^A(v_0), \tilde{//_s}d\tilde B_s\right>_{x_0} |x_t=y\right\},
\hskip 10pt v\in T_{x_0}M,$$
where $\{x_s\}$ is a diffusion on $M$ with generator
 $\half\hbox{trace}\tilde \nabla \hbox{grad -}+L_A$, and $\tilde B$ the
martingale part of the  stochastic  anti-development of
 $\{x_s: 0\le s\le t\}$ using $\tilde \nabla$, a 
Brownian motion  on $T_{x_0}M$, while $\tilde{//_s}$ is parallel
 translation,  and $v_s=\tilde W_s^A(v_0)$ satisfies
$${\tilde D^\prime\over \partial s} v_s =-\half \tilde{\hbox{Ric }^\#} (v_s)
+\tilde \nabla A(v_s)$$
both along the paths of $\{x_s: 0\le s<t\}$ where $\tilde D^\prime$
refers to covariant differentiation using the adjoint connection
 $\tilde \nabla^\prime$.
\end{corollary}

\noindent
{\bf Proof.}
As described in \S 2H we can choose an SDE (\ref{1}) with
 $\breve \nabla =\tilde \nabla$ and then  the generator
is as required by (\ref{118}). If $\xi_\cdot$ is the flow then
by \cite{ELflow}, 
$$d\log p_t(\cdot, y)(v_0)={1\over t}\E \left\{
\int_0^t \left. \left <T\xi_s(v_0), X(x_s)dB_s \right>
\right | \xi_t(x_0)=y\right \}$$
and the result follows  from the theorem and the fact that
$\breve B_t$ given by
$d\breve B_t=\breve{//_t}^{-1} X(x_t)dB_t$ is the martingale
part of the  stochastic  anti-development (defined by the corresponding
 Stratonovich equation).
\hfill//

\bigskip

\noindent{\bf Example 3. }
For the flat left invariant connection on a Lie group $G$ the SDE is
 as described in \S2 G. Then $W_t$ and $T_{x_0}\xi_t$ are equal (there is
no extraneous noise) and they are just right translation by $\xi_t(x_0)$
while $\tilde{//_t}$ is left translation by $\xi_t(x_0)$.

\section{Moment Exponents}

Let $S(t,x)(e)$ be the flow for the vector field $X^e$, and set
$\delta S(t,v)(e)=TS(t,x)(e)(v)$. Let $<,>^\prime$ be a  Riemannian metric
on $M$, not necessarily the induced one from the SDE.  Denote by $|-|^\prime$
the corresponding norm.
Let $\tilde \nabla^\prime$ be a connection compatible with $<,>^\prime$.
 Then

\begin{equation}\label{600}
{d\over dt}|\delta S(t,v)e|^{\prime\,p}=
p\left| \delta S(t,v)e\right|^{\prime \, p-2}
\left<\delta S(t,v)e, {\tilde D^\prime \over dt}\delta S(t,v)e\right>^\prime.
\end{equation}

\noindent
Also

\begin{eqnarray*}
{\tilde D^\prime\over \partial t} \delta S(t,v)(e)
&=&\tilde \nabla^\prime X^e\left(\delta S(t,v)e\right)+
\tilde T^\prime\left( X^e(S(t,x)e), \delta S(t,x)e \right)\\
&=&\tilde \nabla X^e\left(\delta S(t,v)e\right),
\end{eqnarray*}
as for  (\ref{500}) if $\tilde \nabla$ is the adjoint of $\tilde\nabla^\prime$.
Then

\begin{equation}
{d\over dt}|\delta S(t,v)e|^{\prime\,p}=
p\left| \delta S(t,v)e\right|^{\prime\, p-2}
\left<\delta S(t,v)e,
 \tilde \nabla X^e\left(\delta S(t,v)e\right)\right>^\prime.
\end{equation}
At $t=0$,
\begin{equation}\label{30}
{d\over dt}|\delta S(t,v)e|^{\prime\, p}=
p|v|^{\prime\,p-2} <v,\tilde \nabla X^e(v)>^\prime.
\end{equation}
Furthermore
\[\begin{array}{l}
 {d^2\over dt^2}\left. \left|\delta S(t,v)\right|^{\prime\,p} \right|_{t=0}\\
= p(p-2) |v|^{\prime\,p-4}<v,\tilde\nabla X^e(v)>^{\prime\, 2}
+ p|v|^{\prime\, p-2}\left[    |\tilde \nabla X^e(v)|^{\prime\, 2}
+\left<\tilde T\left(\tilde\nabla_v X^e, X^e\right),\, v\right>^\prime \right]\\
+ p|v|^{\prime\,p-2} \left<v, \tilde \nabla^2 X^e (X^e, v)\right>^\prime
+ p|v|^{\prime\,p-2}
 \left<v, \tilde \nabla X^e(\tilde \nabla X^e(v))\right>^\prime.
\end{array}\]

\noindent Set

\[\begin{array}{ll}
H_p(x)(v,v)& = 2<\tilde\nabla A(v),v>^\prime
+\sum_1^m \left<\tilde \nabla^2 X^{i}(X^{i}, v), v\right>^\prime\\
&+\sum_1^m\left<\tilde \nabla X^{i}(\tilde \nabla X^{i}(v)), v\right>^\prime
+ \sum_1^m \left|\tilde \nabla X^{i}(v)\right|^{\prime\, 2}\\
&+\sum_1^m\left[
\left< \tilde T\left(\tilde\nabla_v X^i, X^i\right),\,v\right>^\prime
 +(p-2) {1\over |v|^{\prime\,2}}
\left<\tilde \nabla X^{i}(v), v\right>^{\prime\, 2}\right].\\
\end{array}\]
In terms of the Ricci curvature,
 
\[\begin{array}{ll}
H_p(x)(v,v)&
=  2<\tilde\nabla \left(A+\sum_1^m\tilde \nabla_{X^i} X^i\right)(v),v>^\prime
-<\tilde{\rm Ric}^\#(v),v>^\prime\\
&+\sum_1^m\left[
\left<\tilde T\left(\tilde\nabla_v X^i, X^i\right),\, v\right>^\prime +
  \left|\tilde \nabla X^{i}(v)\right|^{\prime\, 2}
 +(p-2) {1\over |v|^{\prime\, 2}}
\left<\tilde \nabla X^{i}(v), v\right>^{\prime\, 2}\right].\\
\end{array}\]

 From (\ref{30})  and the equation after it we see that
$<v,\tilde \nabla X^i(v)>^\prime$ and $ H_p(x)(v,v)$ are independent of 
the choice of such  connections for fixed $<-,->^\prime$.
 In particular when $<,>^\prime$ is the metric $<,>$ induced by the S.D.E.
the $H_p$ defined here agrees with the one used in 
\cite{application}. 

Taking $\tilde \nabla=\breve \nabla$, we see that
if $\hat \nabla$ is compatible with $<-,->^\prime$,
\begin{eqnarray*}
H_p(x)(v,v)&=&  2<\breve\nabla A(v),v>^\prime
-<\breve{\rm Ric}^\#(v),v>^\prime\\
&& +\sum_1^m \left|\breve \nabla X^{i}(v)\right|^{\prime\, 2}
 +(p-2) \sum_1^m {1\over |v|^{\prime\,2}}
\left<\breve \nabla X^{i}(v), v\right>^{\prime\, 2}.
\end{eqnarray*}

\bigskip

By It\^o's formula (c.f. \cite{ELflour}), we have

\begin{lemma}
Let $\tilde\nabla $ be a connection whose dual connection is metric
for some metric  $<-,->^\prime$. Then for $v_0\in T_{x_0}M$,
\begin{eqnarray*}
|T\xi_t(v_0)|^{\prime \,p}=&&|v_0|^{\prime\, p} +
\int_0^t p|T\xi_s(v_0)|^{\prime\, p-2} <T\xi_s(v_0), 
\tilde \nabla X(T\xi_s(v_0))dB_s>^\prime\\
&&+{p\over 2} \int_0^t |T\xi_s(v_0)|^{\prime\,p-2} 
 H_p(\xi_s(x_0))(T\xi_s(v_0), T\xi_s(v_0))ds,
\end{eqnarray*}
\end{lemma}

\bigskip 

\noindent  Set
\begin{eqnarray*}
h_p(x)&=&\sup_{|v|=1}  H_p(x)(v,v),\\
\underline{h}_p(x)&=&\inf_{|v|=1} H_p(x)(v,v).
\end{eqnarray*}

We can now extend the result proved in \cite{application} for gradient Brownian 
systems:
\begin{proposition}
Suppose $\hat \nabla $ is  metric for some Riemannian metric $<-,->^\prime$.
 Then
\begin{equation}\label{32}
\E \exp^{\half \int_0^t \underline{h}_p(\xi_s(x_0))ds}\le
\E|T_{x_0}\xi_t|^{\prime \,p}
\le n \E \exp^{\half \int_0^t h_p(\xi_s(x_0))ds}.
\end{equation}
\end{proposition}

\noindent {\bf Proof.}
Let $P_N(x): \R^m\to N(x)$ be the orthogonal projection.
 Define
$$A_x: T_xM\oplus N(x)\to T_xM, \hskip 19pt x\in M $$
by 
$$A(u,e)=\breve\nabla X (e)(u).$$
Then $A$ is the shape operator when (\ref{1}) is a gradient system.
For $e\in \R^m$ we have
$$A(u, P_N(x)(e))=\breve \nabla X(e)(u).$$
Note that we can write
$$|T\xi_t(v_0)|^{\prime\,p}=
|v_0|^{\prime\,p} {\varepsilon}(M_t^p)\exp^{a_t^p}$$
for ${\varepsilon}(M_t^p)$ the exponential martingale corresponding
 to $M_t^p$ where
 $$M_t^p= \sum_1^m \int_0^t p\,
 {<\breve \nabla X^i(T\xi_s(v_0)), T\xi_s(v_0)>_{\xi_s(x_0)}^{\prime}
\over |T\xi_s(v_0)|^{\prime\,2}} dB_s^i$$
and for
 $$a_t^p ={p\over 2}\int_0^t 
{H_p(\xi_s(x_0))(T\xi_s(v_0), T\xi_s(v_0))\over
 |T\xi_s(v_0)|^{\prime\,2}}ds.$$
Now we are in the situation of \cite{application}
and the same proof, by the Girsanov transformation as used there,
leads to (\ref{32}).
\hfill  //

\bigskip
It is worth mentioning that since (\ref{30}) and the equation after
it is invariant under choice of connections we see that if $\breve \nabla$
is torsion skew symmetric then
\begin{eqnarray*}
\breve{\hbox{Ric}}_x(v,v)&=&\hbox{Ric}_x(v,v)- \sum_{1}^m |\nabla X^i(v)|^2
+\sum_1^m |\breve \nabla X^i(v)|^2\\
&=& \hbox{Ric}_x(v,v)-\sum_{1}^n |\nabla X^i(v)|^2
\end{eqnarray*}
because for such connections $ \sum_{i=1}^m\nabla X^i(X^i)=0$ 
by Corollary \ref{2.5}  and
$\nabla X^i(v)=\breve\nabla X^i(v)$ for $i>n$ since $X^i(x_0)=0$.
In particular

\begin{proposition}\label{4.3}
The Ricci curvature of any torsion skew symmetric connection
$\tilde \nabla$ is majorized by that of  the corresponding
Levi-Civita connection. Equality holds everywhere if and
only if $\tilde \nabla$ is Levi-Civita.
\end{proposition}

\section{The generator on differential q-forms}

Let $\phi$ be a differential q-form and  $\xi_t^* \phi$
 its pull back  by our flow $\xi_t(-)$.
This gives rise to a semigroup on bounded q-forms \cite{ELflow}:
$$P_t\phi=\E \xi_t^* \phi,$$
 i.e.   if $v=(v_1,\dots, v_q)$ is a q-vector in $\bigoplus^q T_xM$, 
$P_t\phi(v)=\E\phi(T\xi_t(v_1), \dots, T\xi_t(v_q))$.

Its infinitesimal generator $\A^q$ is given by:
$$\A^q\phi=\left({1\over 2}\sum_1^m L_{X^i}L_{X^i}+L_A\right)\phi,$$
where $L_A$ denotes Lie differentiation in the direction of $A$.

\bigskip

Let $i_{A}\phi$ be the interior product of $\phi$ by $A$, which is a q-1 form
defined by: $i_A\phi(v_1, \dots, v_{q-1})=\phi(A, v_1,\dots, v_{q-1})$. Set
\begin{equation}
\bar \delta \phi= -\sum_1^m i_{X^i} \hat\nabla\phi(X^i).
\end{equation}
Then it is easy to see that $\bar \delta \phi= -\sum_1^m i_{X^i}L_{X^i}\phi$
and
$$\sum_1^m L_{X^i}L_{X^i}\phi= -\bar \delta d\phi- d\bar \delta\phi$$
for $d$ the exterior differentiation.

\bigskip

 There is also a Weitzenb\"ock formula:
\begin{equation}
\A^q\phi=\half {\rm trace}\hat \nabla^2\phi -\half \breve R^q(\phi)+
L_A(\phi)
\end{equation}
where $\breve R^q$ is the  zero order operator on q-forms obtained
 algebraically (e.g. via annihilation and creation operators
as in \cite{Cycon-Froese-Kirsch-Simon} or see \cite{ELflour})
 from the curvature tensor $\breve R$ of
$\breve \nabla$ just as the usual Weitzenb\"ock terms are  obtained
from the curvature of the Levi-Civita connection. In particular
for a 1-form $\phi$,
$$\breve R^1(\phi)(v)=\phi(\breve{{\rm Ric}}^\#(v)), \hskip 19pt v\in T_xM.$$
The case of 1-form is straightforward, or can be seen from Theorem 
\ref{th:2C}. For details of the general case and further discussions see
\cite{EL-LJ-LI}.

\section*{Acknowledgment}
This research was supported by SERC Grant GR/H67263 and 
EC grant SC1*-CT92-0784.

\section*{Appendix I \hskip 8pt  The Curvature Tensor}

To calculate the curvature tensor $\breve R$ we will use the expression
in Lemma \ref{le:2.1} (ii) for $\breve \nabla$. Thus if $U, V, W$
are vector fields

\begin{eqnarray*}
\breve \nabla_U \breve \nabla_V W &=&[U, \breve \nabla_V W]+
\sum_1^m [X^i,U]<\breve \nabla_V W, X^i>\\
&=&[U, [V,W] ]+\sum_1^m [U, [X^i, V] ]<W, X^i>\\
&& \hskip 4pt + \sum_1^m\, [X^i, V] \, d<W, X^i>(U(\cdot))
+ \sum_1^m\, [X^i, U] <\breve \nabla_V W, X^i>.
\end{eqnarray*}

Applying Jacobi's identity twice we see

\begin{eqnarray*}
\breve R(U,V)W:&=& \breve \nabla_U \breve \nabla_V W -
\breve \nabla_V \breve \nabla_U W -\breve \nabla_{[U,V]} W\\
&=&  \sum_1^m \left\{\, [X^i, V] \,d<W, X^i>(U(\cdot)) \right.\\
&& \hskip 8pt - \left.   [X^i, U] \,d<W, X^i>(V(\cdot)) \, \right\}\\
&& \hskip 5pt 
+\sum_1^m \left\{  [X^i, U] <\breve \nabla_V W, X^i>
- [X^i, V] <\breve \nabla_U W, X^i> \right\}.
\end{eqnarray*}

Now take $U=Z^u$, $V=Z^v$, $W=Z^w$ for $u, v ,w \in T_{x_0} M$. Then
\begin{eqnarray*}
\breve R(u,v)w&=&
\sum_1^m \left\{ [X^i, Z^v]<\breve \nabla_u X^i, w> 
-[X^i, Z^u] <\breve \nabla_v X^i, w>\right\}\\
&=& \sum_1^m \left\{ -\nabla_v X^i<\breve \nabla_u X^i, w>
+\breve \nabla_u X^i <\breve \nabla_v X^i, w>\right\}
\end{eqnarray*}
since the torsion terms vanishes when summed in conjection with
the terms which involve $\breve \nabla X^i$. Thus

\bigskip

\noindent
{\bf Proposition A1 } {\it  If $u, v, w \in T_{x_0}M$ then

$$\breve R(u,v)(w)
=\sum_{i=1}^m\breve\nabla_u X^i<\breve\nabla_v X^i, w> - 
\sum_{i=1}^m\breve\nabla_v X^i<\breve\nabla_u X^i, w>. $$	}

\noindent
{\bf Corollary A2} 
$$<\breve R(u,v)w, z>
=-\sum_1^m  <\breve \nabla_u X^i \wedge \breve \nabla_v X^i,
 w\wedge z>_{\Lambda^2 T_xM}. $$

\noindent
\noindent{\bf Remark:}
  For $A$ the 'shape operator' defined in \S 4,
the proposition gives
$$\breve R(u,v)w=\hbox{ trace } \left\{
A(u,-) \left<A(v,-), w\right>-  A(v,-)\left<A(u,-),w\right>\right\}$$
which reduces in the gradient case to Gauss's equation for the 
curvature of a submanifold in $\R^m$ (e.g. p. 23 \cite{Kobayashi-NomizuII}).

\section*{Appendix II}

There is a direct correspondence between stochastic flows and
 Gaussian measures $\gamma$ on the space $\Gamma(TM)$ of vector fields on 
$M$, \cite{Baxendale82}, \cite{LeJan-Watanabe82}, \cite{Kunitabook}.
 The latter 
is determined by its reproducing kernel Hilbert space 
(Cameron-Martin space), a Hilbert space $H$ of vector fields on $M$,
 together with its mean $\bar \gamma$, a vector field on $M$
 \cite{Baxendale76}. For the flow
corresponding to our S.D.E. (\ref{1}), the measure $\gamma$ is the image
measure of the standard Gaussian measure on $\R^m$ by the map
\hskip 10pt $e\mapsto X^e$ \hskip 10pt
from $\R^m$ to vector fields on $M$ shifted by $\bar \gamma$,
in this case the vector field $A$. The space $H$ is just
 $\{X^e: e\in \R^m\}$ with quotient inner product. However in general
 $H$ may be infinite dimensional e.g. for isotropic stochastic flows
\cite{LeJan85}.

\bigskip

Nevertheless given such $H$ and vector fields $\bar \gamma$, if $\gamma_0$
is the corresponding centered Gaussian measure and $\{W_t: t\ge 0\}$
the Wiener process on the space of vector fields with $W_1$ distributed
as $\gamma_0$, the corresponding stochastic flow is obtained
as the solution flow of
$$dx_t=\rho_{x_t}\circ dW_t+\bar \gamma(x_t)dt$$
where $\rho_x: \Gamma(TM) \to T_xM$ is the evaluation map
(assuming sufficient regularity), see \cite{ELflow}. This reduces
the situation to that discussed above with $\R^m$ replaced by
the possibly infinite dimensional Hilbert space $H$. Assume
non-degeneracy, so $\rho_x$ is surjective for each $x$, and
let $M$ have the induced Riemannian metric. It is worth
noting that the adjoint $Y(x): T_xM \to \R^m$ of $X(x)$ is now
replaced by the adjoint of $\rho_x: H\to T_xM$ which is essentially
the reproducing kernel of $H$, i.e. the covariance of $\gamma$:
$$\rho_x^*(v)=k(x,\cdot)(v) \in H$$
where 
\begin{equation}
<k(x,\cdot)v, h>_H= <h(x), v>_x, \hskip 18pt x\in M, v\in T_xM.
\end{equation}
In particular 
$$Z^v=k(x_0, \cdot) v, \hskip 18pt v\in T_{x_0}M.$$
and for a vector field $Z$ on $M$ our basic definition (\ref{5}) becomes
\begin{equation}\label{100}
\breve \nabla Z(v)=d[ k(\cdot, x_0)Z(\cdot)](v)
\end{equation}
treating $y\mapsto k(y,x_0)Z(y)$ as a map from $M$
to $T_{x_0}M$. The defining condition (\ref{4}) for $\breve \nabla$
can be written $\breve \nabla (k(x_0,\cdot)v)w=0$ all $v, w\in T_{x_0}M$
 all $x_0\in M$. 

\bigskip

In terms of expectation with respect to our basic
Gaussian measure $\gamma_0$, treating vector fields $W$ as a random field, 
equation (\ref{10}) for $\breve \nabla Z$ becomes
\begin{equation}
\breve \nabla Z(v)
={d\over dt}\left. \E W(x_0) <Z(\sigma(t), W(\sigma(t))>_{\sigma(t)}
\right|_{t=0}
\end{equation}
$$k(x,\cdot)v=\E <W(x), v>_x W(\cdot)$$
and in terms of conditional expectations
$$k(x, y)W(x)=\E\left\{W(y)|W(x)\right\}  \hskip 5pt \in T_yM$$
$$k(x, \cdot)v=\E \left\{ W|W(x)=v\right\},$$
giving 
$$\breve\nabla Z(v)={d\over dt}\E
\left\{ W(x_0) \left|\, W(\sigma(t))=Z(\sigma(t))\right\} \right|_{t=0}.$$

\end{document}